\documentclass[12pt,leqno]{article}
\usepackage{amsmath,amsfonts,amssymb,graphicx,epsfig}
\usepackage[pdfpagemode=None,colorlinks=true,linkcolor=blue,citecolor=blue]{hyperref}

\newtheorem{rema}{Remark}

\newtheorem{lemm}{Lemma}
\newtheorem{theo}{Theorem}
\newtheorem{coro}{Corollary}
\newtheorem{exem}{Example}
\newenvironment{proof}[1]{\noindent\emph{Proof.} #1}{\hfill$\square$}

\newcommand{\Z}[1][]{\ensuremath{{\mathbb{Z}^{#1}} }}
\newcommand{\C}[1][]{\ensuremath{{\mathbb{C}^{#1}} }}
\newcommand{\R}[1][]{\ensuremath{{\mathbb{R}^{#1}} }}

\renewcommand{\S}[1][]{\ensuremath{{\mathbb{S}^{#1}} }}

\def\Im{ \mathrm{Im}\, }
\newcommand{\<}{\langle}
\renewcommand{\>}{\rangle}
\newcommand{\ga}{\gamma}
\newcommand{\pa}{\partial}
\newcommand{\si}{\sigma}
\newcommand{\al}{\alpha}
\date{}
\author{Henri Anciaux, Ildefonso Castro and Pascal Romon}

\title{Lagrangian submanifolds foliated by
$(n-1)$-spheres in $\R^{2n}$}

\begin{document}
 \maketitle

\begin{abstract}
We study Lagrangian submanifolds foliated by $(n-1)$-spheres in 
$\R^{2n}$ for $n \geq 3$. We give a parametrization valid for such submanifolds,
and refine that description when the submanifold is 
special Lagrangian, self-similar
or Hamiltonian stationary. In all these cases, the submanifold is 
\emph{centered}, i.e. invariant under the action of $SO(n)$. It suffices then
to solve a simple ODE in two variables to describe the geometry of the solutions. \\


\end{abstract}

 \section*{Introduction}
 Lagrangian submanifolds constitute a distinguished subclass in
 the set of $n$-dimen\-sion\-al
 submanifolds of $\R^{2n}:$ a submanifold is \em Lagrangian \em
if it has dimension $n$ and is isotropic with respect to the symplectic
form $\omega$ of $\R^{2n}$ identified with $\C^n$. Such submanifolds
  are locally
 characterized as being graphs of the gradient of a real map on
 a domain of $\R^{n}$, however to give a general classification
 of them is not an easy task.

 In the case of surfaces i.e. $n=2$, powerful techniques such as
  Weierstrass representation formulas may be used
 to gain some insight into this question (\emph{cf} \cite{HR}).
In higher dimension, some work has been done in the direction of characterizing
Lagrangian submanifolds with various intrinsic and extrinsic geometric assumptions,
see \cite{Ch} for a general overview and more references.


 \bigskip

Recently, D. Blair has studied in \cite{Bl}
Lagrangian submanifolds of $\C^n$ which are 
foliated by $(n-1)$-planes.
 In the present paper, we are devoted to study those Lagrangian submanifolds of
 $\R^{2n}$, 
 which are foliated by Euclidean $(n-1)$-spheres.
 In the following, we shall for brevity denote them by
  \em $\sigma$-submanifolds. \em
 We first observe that
  any isotropic round $(n-1)$-sphere spans a $n$-dimensional 
Lagrangian subspace, 
 a condition which breaks down
 in dimension
 $2$. Hence we restrict our attention to the case $ n \geq 3$
and we give a characterization of $\si$-submanifolds
 as images of an immersion of a particular form involving
 as data a planar curve and a $\R^n$-valued curve (see Theorem 1).

 Next, we focus our attention on several curvatures equations
 and characterize those among which are $\si$-submanifolds.
  The most classical of these equations is the \em minimal
  submanifold equation, \em involving the mean curvature vector 
$\vec{H}$
  $$\vec{H}=0.$$
  A Lagrangian submanifold
  which is also minimal satisfies in addition a very striking
  property: it is calibrated (\emph{cf} \cite{HL}) and therefore
  a minimizer; these submanifolds are called
  \em special Lagrangian. \em Many special Lagrangian submanifolds
  with homogeneity properties
  have been described in \cite{CU2,Jo}.
 In this context we recover
 the characterization of the Lagrangian catenoid as the only (non flat)
 special Lagrangian $\sigma$-submanifold (see \cite{CU1}).

  We shall call \em self-similar \em a submanifold satisfying
   $$ \vec{H}+ \lambda X^{\perp} =0,  $$
where $X^{\perp}$ stands for the projection of the position vector
$X$ of the submanifold 
on its normal space and $\lambda$ is some real constant (\emph{cf}  \cite{An2}). 
Such a submanifold has the property that its
evolution under the mean curvature flow is a homothecy (shrinking
to a point if $\lambda > 0$ and expanding to infinity if $\lambda
< 0$). Here we shall show that there are no more self-similar
$\sigma$-submanifolds than the ones described in \cite{An2}.

The third curvature equation we shall be interested in deals with
 the \em Lagrangian angle \em $\beta$
 (\emph{cf} \cite{Wo} for a definition),
 a $\R / 2 \pi \Z$-valued
 function which is defined up to an additive constant on any
 Lagrangian submanifold and satisfies 
$J \nabla \beta=n \vec{H}$, where $J$ is the complex structure in 
$\C^n$ and $\nabla$ is the gradient of the induced metric on the
submanifold. Following \cite{Oh}, 
 we call \em Hamiltonian stationary \em a Lagrangian submanifold
 which is critical for the volume functional under Hamiltonian
 deformations (generated by vector fields $V$ such that 
$V \lrcorner \; \omega$ is exact).
 It turns out that the corresponding Euler-Lagrange equation is
 $$\Delta \beta = 0.$$
 In others words a Lagrangian submanifold
 is Hamiltonian stationary if and only if its Lagrangian angle
 function is harmonic (for the induced metric).
 Many Hamiltonian stationary surfaces in $\C^2$ have been described
 in \cite{HR} and \cite{An1}, but in higher dimensions, very few examples
 were known so far (for example  the Cartesian product
 of round circles
 $a_1 \S^1 \times \ldots \times a_n \S^1$). In this paper we shall 
describe all Hamiltonian stationary $\si$-submanifolds.

 \bigskip

 The paper is organized as follows: the first section gives the
 proof of characterization of $\si$-submanifolds (Theorem 1).
 In the second one, we
 compute for this kind of submanifold the Lagrangian angle and
 the mean curvature vector.
 In Section 3, we obtain as a corollary the fact that special
 Lagrangian or self-similar $\si$-submanifolds are
 centered (i.e. invariant under the standard 
action of the special orthogonal group $SO(n)$, see Example 1 in Section 1
 for a precise definition),
 thus the only examples of such submanifolds are those described
 in \cite{CU1} (for special Lagrangian ones) and in \cite{An2} (for self-similar
 ones).
 In Section 4, we show that the only Hamiltonian stationary
 $\si$-submanifolds are also centered (Theorem \ref{HS}),
 and in Section 5 we describe
 them in details (Corollary 1).

\bigskip

\em Acknowledgments: \em The authors would like to thank
Francisco Urbano for his numerous valuable comments and suggestions
and
Bruno Fabre for the proof of Lemma \ref{Alg}.

 \section{Characterization of Lagrangian submanifolds foliated
by $(n-1)$-spheres}
 On $\C^n \simeq \R^{2n}$, $n \geq 3$, with coordinates
 $\{z_{j}=x_{j}+iy_{j}, \, 1 \leq j \leq n \}$ equipped with
 the standard Hermitian form $\<.,.\>_{\C^n}$
and its associated symplectic form
  $ \omega:=\sum_{j=1}^{n} dx_{j} \wedge dy_{j}=\Im\<.,.\>_{\C^n}$,
 we consider submanifolds of dimension $n$ foliated by round
 spheres $\S^{n-1}$. Locally, they can be parameterized by immersions

  $$ \begin{array}{lccc} \ell : & I \times \S^{n-1}  &\to& \C^n \simeq \R^{2n}
 \\
&  (s , x) & \mapsto & r(s)M(s)x + V(s) \end{array}$$
where
\begin{itemize}

    \item  $I$ is some interval,
\item $\forall s \in I$, $ r(s) \in \R^{+}, M(s) \in SO(2n)$, $V(s) \in \C^n$,

\item $x \in \S^{n-1} \subset \R^n \subset \C^n$, where we note
  $\R^{n}=\{ y_{j}=0 , \,  1 \leq j \leq n \}$.

\end{itemize}

\begin{lemm}
  A necessary condition for the immersion $\ell$ to be Lagrangian is
  that for any $s$ in $I$, the (affine) $n$-plane containing the leaf
  $\ell(s,\S^{n-1})$
  is Lagrangian; in
  other words, $M \in U(n)$, where $U(n)$ is embedded as a subgroup
of $SO(2n)$.
\end{lemm}
\begin{proof}
For a fixed $s \in I$, the leaf $\ell(s,\S^{n-1})$ spans exactly the affine
$n$-plane  $V+M\R^n$. So we
have to show that for two independent vectors $W, W' \in \R^n$,
$\omega(MW,MW')=0$.
Since $n \geq 3$, there exists $x \in \S^{n-1}$ such that
$W$ and $W'$ are tangent to $\S^{n-1}$ at $x$. Thus
$rMW$ and $rMW'$ are tangent to the leaf $\ell(s,\S^{n-1})$
at $\ell(s,x)$. The tangent space to this leaf is isotropic (being
included in the tangent space to $\ell(I \times \S^{n-1})$).
Therefore  $\omega(MW,MW')=0$. Finally, $M$ is an isometry mapping 
a Lagrangian $n$-plane to another 
Lagrangian $n$-plane, hence a unitary transformation.
\end{proof}

\begin{rema} This crucial lemma does not hold in dimension two: we can
    produce examples of Lagrangian surfaces foliated by circles such
    that the planes containing the circular leaves are not Lagrangian
    (cf Example 3).
\end{rema}

Now, we know that for any $\xi \in \R^n$ orthogonal to
$x$, $M\xi$ is tangent to the submanifold at $\ell(s,x)$
(because it is tangent to the leaf) so
we have:
$$\omega \big(\ell_s, M \xi \big)=0,$$
where
$\ell_s=\pa \ell / \pa s = \dot{r}M x + r \dot{M}x+\dot{V}$.
 Along the paper the dot ` $\dot{}$ '  will be a shorthand
 notation for the derivative with respect to $s$. 
 We compute:
$$  \omega(\dot{r} M x, M \xi )
      + \omega( r \dot{M} x, M \xi)+\omega(\dot{V}, M \xi)=0$$
$$ \Leftrightarrow
   \Im\< \dot{r} Mx , M \xi \>_{\C^n}
      + \Im \< r \dot{M}x,  M \xi \>_{\C^n}+\Im\<\dot{V}, M \xi \>_{\C^n}=0$$
$$ \Leftrightarrow
     \dot{r} \Im\< x , \xi \>_{\C^n}
    + r \Im \<M^{-1} \dot{M}x, \xi \>_{\C^n}+\Im\<M^{-1}\dot{V}, \xi \>_{\C^n}=0$$
$$ \Leftrightarrow   r \< \Im(M^{-1} \dot{M}) x, \xi \>_{\C^n}+
      \<\Im(M^{-1}\dot{V}),  \xi \>_{\C^n} =0,$$
because $\xi$ is real.

In the following, we shall note
$b:=\Im (M^{-1} \dot{V}) \in \R^{n}$ and $B:=\Im (M^{-1} \dot{M})$.
The complex matrix $ M^{-1} \dot{M}$ is
skew-Hermitian so $B$ is symmetric.
Then we have the following equation:
\[
\<b,\xi \>_{\R^n} + r \<B x, \xi \>_{\R^n}=0,
  \quad \mbox{for all } (x,\xi)  \in \S^{n-1} \times \R^n
   \mbox{ such that }
     \<x, \xi \>_{\R^n}=0 \tag{$*$}
\]

\bigskip

From now on, we denote the real scalar product
$\<., . \>_{\R^n}$ by
$\<.,.\>$.
Next we have three steps:
\begin{description}
\item[Step 1:] \emph{$b$ vanishes or is an eigenvector for $B$.}\\
If $b \neq 0$, take $x$ collinear to $b$, so that $x=\lambda b$ 
for some real non zero
constant $\lambda$.  From $(*)$, we have
$$\forall \xi \in b^\perp, \< B \lambda b, \xi \> =0.$$
Therefore $Bb$ is collinear to $b$, so there exists $\mu$ such
that $Bb=\mu b$.

\item[Step 2:] \emph{$b$ vanishes.}\\
Suppose $b \neq 0$, then set $\xi = b$ and $x$ orthogonal to $b$,
so $\<b, x \>=0$ and $\< Bb, x \>=0$ (using step~1).
Using the fact that $B$ is symmetric, we write:
$$ \<b, B x \> =0,$$
so
$$\< \xi , Bx \>=0.$$
From $(*)$ we deduce that $\<b, \xi \>=||b||^2=0$.
So $b$ vanishes.

\item[Step 3:] \emph{$B$ is a homothecy.}\\
Now, $(*)$ becomes
$$ \<Bx,\xi \> =0,
  \quad \forall x, \xi \in \S^{n-1} \mbox{ such that }
     \<x, \xi \>=0.$$
This implies that any vector $x$ is an eigenvector for $B$, so $B$ is a
homothecy.
\end{description}
We deduce that $M^{-1} \dot{V} \in \R^{n}$, and that
$M^{-1} \dot{M} \in
\mathfrak{so}(n) \oplus i \R \mathit{Id}$. This implies that $ M \in
SO(n).U(1)$. In other words, we may write $M=e^{i \phi} N, $
where $N \in SO(n)$. Moreover, we observe that $N$ can be fixed to
be $\mathit{Id}$: this does not change the image of the immersion. Then
$M=e^{i \phi} \mathit{Id}$, so $\dot{V}= e^{i \phi} W$, where $W \in \R^{n}$.

\bigskip

Finally, we have shown:

\begin{theo}
    Any Lagrangian submanifold of $\R^{2n}, n \geq 3$,
which is foliated by round $(n-1)$-spheres
is locally the image of an
    immersion of the form:
  $$ \begin{array}{lccc} \ell : 
 & I \times \S^{n-1}  & \longrightarrow  & \C^n \simeq \R^{2n}  \\
&  (s , x) & \longmapsto & 
 r(s)e^{i\phi(s)}x + \int_{s_0}^s e^{i \phi(t)} W(t)dt,
    \end{array}$$
where
 \begin{itemize}
   \item $s \mapsto \ga(s):=r(s) e^{i \phi(s)} $ is a planar curve,
   \item $W$ is a curve from $I$ into $\R^n$,
    \item $s_0 \in I$.
  \end{itemize}
\end{theo}

Geometrically, $r(s)$ is the radius of the spherical leaf and
$V(s):=\int_{s_0}^s e^{i \phi(t)} W(t)dt$ its center. The fiber
lies in the Lagrangian plane $e^{i \phi(s)} \R^n$.

\begin{exem} \em 
If $V=W=0$, the center of each leaf is fixed. In the following, we shall
simply call these submanifolds \em centered. \em
In this case, the
submanifold is $SO(n)$-equivariant (for the following action $z
\mapsto Az, A \in SO(n)$ where $SO(n)$ is seen as a subgroup of
$U(n)$, itself a subgroup of $SO(2n))$.
\end{exem}

\begin{exem}
\em  Assume the curve $\ga$ is a straight line passing through the
origin. Then up to reparametrization 
we can write $\ga(s) =  s e^{i \phi_0}$, where
$\phi_0$ is some constant. This implies that
$$ \ell(x,s)= e^{i \phi_0} \left[ s x + \int_{s_0}^s  W(t)dt \right],$$
 thus the immersion is totally geodesic
 and the image is simply an open subset of the Lagrangian subspace
 $e^{i \phi_0} \R^n$.  \em
\end{exem}

\begin{exem}
\em Assume the centers of the leaves lie on some straight line.
Then there exists some function $u(s)$ such that $V(s)=u(s)a+c$
with $a,c  \in \C^n$. Differentiating, we obtain
$$
\dot{u}(s)a=e^{i \phi(s)}W(s).
$$
As $a$ is constant and $W$ real, it implies that $\phi$ is
constant, so we have $a=e^{i \phi}b$, where $b \in \R^n$.
Then
$$\ell(x,s)=e^{i \phi}(r(s)x+u(s)b)+c$$
so again the immersion is totally geodesic. We notice that this
situation is in contrast with the case of dimension 2, where there
 exists a Lagrangian flat cylinder, which is foliated by
round circles whose centers lie on a line. \em
\end{exem}

\begin{exem}
\em Epicycloids: Assume the centers of the leaves lie on a circle
contained in a complex line. Then there exists some function
$u(s)$ such that $V(s)=e^{i u(s)}a+c$ with $ a,c \in \C^n$.
Differentiating, we obtain
$$W(s)=i \dot{u}(s)e^{i (u(s)-\phi(s))}a.$$
As in the previous example, it implies that $\phi-u$ is constant;
without loss of generality, we can take $u=\phi$ and
 $a=-ib$, where $b \in \R^n$.
Then
$$\ell(x,s)=e^{i \phi(s)}(r(s)x-ib)+c.
$$ \em
\end{exem}


\section{Computation of the Lagrangian angle and of the mean curvature
vector}
\subsection{The Lagrangian angle}
We may assume, and we shall do so from now on, that $\gamma$ is
parameterized by arclength, so there exists $\theta$ such that
$\dot{\gamma}=e^{i \theta}$, and the curvature of $\gamma$ is
$k=\dot{\theta}$. We also introduce, as in \cite{An2},
$\alpha:=\theta-\phi$, so that $\dot{r}=\cos \alpha$ and
 $\dot{\phi}=\frac{\sin \alpha}{r}$.

 Let $(v_2, ... ,v_n)$ be
 an orthonormal basis of $T_x \S^{n-1}$.
Here and in the following, the indices $j$, $k$ and $l$ (and
the sums) run from $2$ to $n$, unless specified.
We have:

$$\ell_s=e^{i \theta} x + e^{i \phi}W \mbox{ and }
\ell_*v_j=r e^{i \phi} v_j.$$

The induced metric $g$ on $I \times \S^{n-1}$ has the following
components with respect to the basis $(\pa_s,v_2, ... ,v_n)$:
$$ g_{11}=1+|W|^2+ 2 \cos \alpha \<W,x\>,$$
$$ g_{1j}=g_{j1}=r\<W,v_j\>,$$
$$ g_{jk}=g_{kj}= \delta_{jk} r^2.$$
We consider the orthonormal basis (for $g$)
$(e_1, ... , e_n)$ defined as follows:
$$e_j:=\frac{v_j}{r}, \quad e_1:= A \partial_s + \sum B_j v_j,$$
where the real numbers $A$ and $B_j$ are uniquely determined
by the orthonormality condition:
\[
A=\left(1+ 2  \cos \alpha \<W,x\>+\<W,x\>^2 \right)^{-1/2}
\; \textrm{ and } \;
B_j= - \frac{A \<W, v_j \>}{r}\;,
\]
so
\[
e_1=  A\left( \pa_s - \sum \frac{ \<W, v_j \> }{r}v_j \right).
\]

\bigskip

This yields an orthonormal basis of the tangent space to the submanifold
 and we may compute the Lagrangian
angle by means of the following formula: 
$e^{i  \beta}=\det_{\C}(\ell_*e_1, ... , \ell_*e_n)$.
Indeed we have
$$\ell_*e_j=e^{i \phi}v_j,$$
and
\begin{eqnarray*}
\ell_*e_1 &=& A(\ell_s-\sum \frac{\<W,v_j\>}{r}\ell_*v_j)
=A(e^{i \theta} x + e^{i \phi}W -\sum \<W,v_j\>e^{i \phi}v_j) 
\\
&=& A(e^{i \theta}+ e^{i \phi}\<W,x\>)x .
\end{eqnarray*}

From the above remarks we deduce that:
\begin{eqnarray*}
\beta &=& \mathrm{Arg}(e^{i \theta}+ e^{i \phi}\<W,x\>)+(n-1)\phi
\\
&=& \mathrm{Arg}(e^{i \alpha}+ \<W,x\>)+n\phi .
\end{eqnarray*}


\bigskip
In the centered case $W=0$ (Example 1),
we find that $\beta=n \phi +\alpha$. This is the only
case where the Lagrangian angle is constant on the leaves.
 We recall that a necessary and sufficient condition for a
Lagrangian submanifold to be a special Lagrangian one is that the
Lagrangian angle be locally constant. So our computation shows
that a special Lagrangian $\si$-submanifold must be centered.
Moreover, it has been shown in \cite{An2} that the only such
submanifolds are pieces of the Lagrangian catenoid (\emph{cf}
\cite{CU1} for a complete description), so we recover one of the
theorems of \cite{CU1}:

\bigskip
\noindent \textbf{Theorem}
  \em A (non flat) special Lagrangian submanifold which is foliated 
by $(n-1)$-spheres is congruent to a piece
of the Lagrangian catenoid. \em

\subsection{Computation of the mean curvature vector}
We first calculate the second derivatives of the immersion.
$$ \ell_{ss}=ik e^{i \theta} x +
           \left( \frac{i \sin \alpha}{r}W+\dot{W} \right)e^{i \phi},$$
$$ v_j(\ell_s)=e^{i \theta} v_j,$$
$$ \ell_{v_j v_k}=-\delta_{jk}  r e^{i \phi} x.$$
Then we obtain the following expressions:
$$\<\ell_{ss},J\ell_s \>= k+k \cos \alpha \<W,x\> +
  \frac{\sin \alpha \cos \alpha}{r} \<W,x\> +
  \sin \alpha \< \dot{W},x\> + \frac{\sin \alpha}{r} |W|^2,$$
$$\<\ell_{ss},J\ell_*v_j \>=\<v_j(\ell_s),J\ell_s\>=\sin \alpha \<W,v_j\>,$$
$$ \<\ell_{v_j v_j},J\ell_s\>=\<v_j(\ell_s),J\ell_*v_j\>=r \sin \alpha,$$
$$ \<\ell_{v_j v_k},J\ell_*v_k\>=0.$$
Thus we have, using the property that the tensor
$C:=\<h(.,.),J.\>$ is
totally symmetric for any Lagrangian immersion:
\begin{eqnarray*}
C_{111} &=& \< h(e_1,e_1), J\ell_*e_1 \> 
\\
&=& A^3 \left\< h \left(\pa_s - \frac{1}{r} \sum \<W,v_j \> v_j,
   \pa_s - \frac{1}{r} \sum \<W,v_j \> v_j \right),
   J\ell_s - \frac{\sum \<W, v_j \> J\ell_*v_j }{r} \right\>
\\
&=& A^3\left(k+k \cos \alpha \<W,x\> +
  \frac{\sin \alpha \cos \alpha}{r} \<W,x\> +
  \sin \alpha \< \dot{W},x\> \right.
\\
&&+ \left. \frac{\sin \alpha}{r} |W|^2
+\frac{\sin \alpha}{r} \sum \<W,v_j\>^2
   - 2 \frac{ \sin \alpha }{r}\sum \<W,v_j\>^2 \right)
\\
&=& A^3 \left(k+k \cos \alpha \<W,x\> +
  \frac{\sin \alpha \cos \alpha}{r} \<W,x\> +
  \sin \alpha \< \dot{W},x\> +
 \frac{\sin \alpha}{r}\<W,x\>^2 \right),
\end{eqnarray*}
$$ C_{jj1}=\< h(e_j,e_j), J\ell_*e_1 \> = \frac{A}{r^2}
 \left\<h(v_j,v_j),J\ell_s-\frac{1}{r}\sum \<W,v_k\>J\ell_*v_k \right\>=
\frac{A \sin \alpha}{r}, $$
\begin{eqnarray*}
C_{11j} &=& \< h(e_1,e_1), J\ell_*e_j \> =
  \frac{A^2}{r} \left(\sin \alpha \<W,v_j \>
-\frac{2}{r} \sum_k \<W,v_k\> r \sin \alpha \delta_{jk} \right)
\\
&=& - \frac{A^2 \sin \alpha \<W, v_j\> }{r},
\end{eqnarray*}
$$ C_{jjk}=\<h(e_j,e_j),J\ell_*e_k  \>= 0.$$
So finally we have the following:
\begin{eqnarray*}
n\vec{H} &=& A \left[A^2 \left(k+k \cos \alpha \<W,x\> +
  \frac{\sin \alpha \cos \alpha}{r} \<W,x\> \right. \right.
\\
&&+ \left. \left. \sin \alpha \< \dot{W},x\> +
 \frac{\sin \alpha}{r}\<W,x\>^2 \right) +
     \frac{(n-1) \sin \alpha}{r}\right] J\ell_*e_1
 \\
 && - \sum  \frac{A^2 \sin \alpha \<W, v_j\> }{r} J\ell_*e_j.
 \end{eqnarray*}

In Section 5, we shall use the following notations:
$ nJ \vec{H}= a e_{1} + \sum a_{j} e_{j}$, where
\begin{multline}
a = -A^{3} \left( k+ \left( k \cos \alpha  +
  \frac{\sin \alpha \cos \alpha}{r} \right) \<W,x\>
    + \sin \alpha \< \dot{W},x \>
 +\frac{\sin \alpha}{r} \<W,x\>^2  \right) 
\\
- (n-1)\frac{A \sin \alpha}{r}
\tag{$H_0$}
\end{multline}
\begin{equation}
a_{j}=  A^2 \frac{\sin \alpha \<W, v_j\>}{r}. \tag{$H_{j}$}
\end{equation}


\section{Application to self-similar equations}
\begin{theo}
 In the class of Lagrangian submanifolds
which are foliated by $(n-1)$-spheres, there
 are no self-translators and the only self-shrinkers/expanders are the
 centered ones described in \cite{An2}.
\end{theo}
\begin{proof}
The self-translating equation is $ \vec{H}=V^{\perp}$
for some fixed  vector $V \in \C^n$.
In particular,
$$  \< \vec{H}, J\ell_*e_j \>= \<V, J\ell_*e_j \>$$
$$ \Longleftrightarrow - \frac{A^2 \sin \alpha}{nr} \<W,v_j\>=
     \<V, J r e^{i \phi}e_j\>$$
$$ \Longleftrightarrow - \frac{A^2 \sin \alpha}{nr} \<W,v_j\>=
       \<V, Je^{i \phi}v_j\>$$
$$ \Longleftrightarrow - \frac{\sin \alpha}{nr}\<W,v_j\>=
    \left(1+2 \cos \alpha \<W,x\>+ \<W,x\>^2 \right)
    \<V, Je^{i \phi}v_j\>.$$
Differentiating this last equation with respect to $v_k$,  $k \neq
j$ (this is possible since $n \geq 3$), we obtain
$$0=\Big(1+2\<W,v_k\>(\cos \alpha+\<W,x\>)\Big)
 \<V, Je^{i \phi}v_j \>.$$
We now observe that the set of points $(x,v_k)$ for which the
first factor in the r.h.s. term of the above equation
vanishes is of codimension $1$ at
most in the unit sphere bundle over $\S^{n-1}$ and apart from this
set, $\<V, Je^{i \phi}v_j \>$ vanishes.
 Coming back to the third
equality of the above equivalences,
 this yields that either $\sin \alpha$ vanishes or
$\<W,v_j\>$ does. In the first case the
 curve $\gamma$ is a line passing through the origin.
 Then we know by Example 2 that the immersion is totally geodesic
 so the image is a Lagrangian subspace. In
 the other case, we know that for $j \neq k$,
 $\<W,v_j\>$ vanishes on some dense open subset of
 points $(x,v_k)$ in the unit sphere bundle over $\S^{n-1}$.
This shows that $W$ vanishes identically. Then it is clear that so
does $V^{\perp}$. This implies that the only solution of the
equation is for vanishing $\vec{H}$, which is the trivial, minimal
case.

\bigskip
The self-shrinking/expanding equation is
$ \vec{H}+ \lambda X^{\perp}=0$,
where $\lambda$ is some real constant.
This implies
$$\< \vec{H}, J\ell_*e_j\> + \lambda \<\ell, J\ell_*e_j\>=0$$
$$ \Leftrightarrow -\frac{A^2 \sin \alpha}{nr} \<W,v_j\> +
   \lambda \left\< \int e^{i \phi}W, Je^{i \phi} v_j \right\> =0$$
$$\Leftrightarrow \frac{A^2\sin \alpha}{nr}\<W,v_j\> =  \lambda
     \left\< \Im( e^{-i \phi} \int e^{i \phi}W), v_j \right\>.$$

 The quantity $\Im(e^{-i \phi} \int e^{i \phi}W)$ depends only
on $s$, so the same argument as above holds, and we deduce that
either $\gamma$ is a line passing through the origin (totally
geodesic case) or $W$ vanishes, which is the centered case
treated in \cite{An2}.
\end{proof}

\section{Hamiltonian stationary $\si$-submanifolds}

The purpose of this section is to prove the following
\begin{theo} \label{HS}
 Any Hamiltonian stationary Lagrangian submanifold foliated
by \hbox{$(n-1)$}-spheres must be centered.
\end{theo}
\begin{proof}
Let $\ell$ be a parametrization of such a submanifold.
We follow the same notations than in Section 2.
 We are going to show that if  $\Delta \beta$ vanishes, then either
$\ga$ is a straight line, so as we have seen in Example 2, we are
in the totally geodesic case, which is in particular centered,
or $W$ vanishes (centered case). The proof is based on an
analysis of the quantity $f(x):=A^{-6} \Delta \beta$
 which turns to be polynomial. Its
  expression is given in the next lemma and the computation
is detailed in Appendix.

\begin{lemm} \label{computation}
 For any
fixed $s$, $f:= A^{-6} \Delta \beta$
 is a polynomial 
  in the three variables
  $ \<W,x\>, \<\dot{W},x\>$ and $\<\ddot{W},x\>$.
Indeed we have: $f=(I)+(II)+(III)+(IV)$ with
\begin{eqnarray*}
(I)&:=&  
3 B  \Big(\cos \alpha \< \dot{W},x\>
    - \dot{\alpha} \sin \alpha \<W,x\>
    + \<W,x\> \< \dot{W},x\> \Big)
\\
&&-
    A^{-2} \left[\dot{B} - (n-1) \frac{ \sin \alpha}{r}
  \Big(\cos \alpha \< \dot{W},x\>
    - \dot{\alpha} \sin \alpha \<W,x\>
    + \<W,x\> \< \dot{W},x\> \Big) \right]
\\
&&-A^{-4} (n-1) \pa_{s} \left(\frac{\sin \alpha}{r} \right),
\\
(II) &:=& 
-3  B \frac{ \cos \alpha +\<W,x\>}{r}(|W|^2-\<W,x\>^2)
\\
&&+\frac{A^{-2}}{r} \Bigg[\left(  k \cos \alpha
  - (n-2) \frac{\sin\alpha\cos\alpha}{r}
  - (n-3) \frac{ \sin \alpha}{r}  \<W,x\>
  \right)  (|W|^2-\<W,x\>^2) 
\\
&&+ \sin \alpha \Big(\<\dot{W},W\>-\<\dot{W},x\> \<W,x \> \Big)
  \Bigg],
\\
(III) &:=& 
- A^{-2} \frac{\sin \alpha}{r^2}
     (\cos \alpha + \<W,x\> ) \Big( |W|^{2}- \<W,x\>^{2} \Big) -
     A^{-4}(n-1) \frac{ \sin \alpha}{r^2} \<W,x\>,
\\
(IV) &:=& 
-A^{-2}B \frac{n-1}{r} (\cos \alpha + \<W,x\>)-
   A^{-4} \frac{(n-1)^{2} \sin \alpha}{r^2}  (\cos \alpha + \<W,x\>)
\end{eqnarray*}
and
$$ B:=k+ \left(k \cos \alpha  +
  \frac{\sin \alpha \cos \alpha}{r} \right) \<W,x\>+
  \sin \alpha \< \dot{W},x\> +
 \frac{\sin \alpha}{r}\<W,x\>^2 .$$
\end{lemm}

\begin{proof}
See Appendix.
\end{proof}

\begin{lemm} \label{Wseconde}
   $$ f(x)= -\Big( \sin \alpha +
     2 \sin \alpha \cos \alpha  \<W,x\>
     + \sin \alpha  \<W,x\>^{2} \Big) \<\ddot{W},x\> +
     R \Big(\<W,x\>,\<\dot{W},x\>\Big),$$
    where $R(.,.)$ is a polynomial in its two variables.
 \end{lemm}

\begin{proof}
In the previous  expression of $f$, we
see that there are no contributions to the terms in
$\<\ddot{W},x\>$ apart from the $\dot{B}$ term in $(I)$. Then we
compute
$$\dot{B}= \sin \alpha \<\ddot{W},x\> + S\Big(\<W,x\>,\<\dot{W},x\>\Big),$$
where $S(.,.)$ is a polynomial in its two variables.
We deduce that the only term of $\<\ddot{W},x\>$ in $f$ is the
following:
$$ - A^{-2} \sin \alpha \<\ddot{W},x\>,$$
so we have our claim.
\end{proof}

\begin{lemm} \label{ordre5}
The polynomial $f$ has total degree at most $5$, and 
in that case its leading term is
$\frac{-n^2+n+2}{r^2} \sin \alpha \<W,x\>^5$. 
\end{lemm}
\begin{proof}
Remembering that $A^{-2}$, $B$ and
$\dot{B}$ are polynomials of
 degree at most $2$ in $\<W,x\>$ and $\<\dot{W},x\>$,
 we first observe that there are no terms
 of order $5$ in  $(I)$.  Then we check
that in $(II)$, $(III)$ and $(IV)$ there are no terms in
$\<\ddot{W},x\> \<W,x\>^4$ (this has already been observed in the
previous lemma) or in $\<\dot{W},x\> \<W,x\>^4$. Finally, the
  coefficient of $\<W,x\>^2$ in $A^{-2}$ and $B$ are respectively
  $1$ and $ \frac{\sin \alpha}{r}$,
 so summing up, we obtain that the coefficients of  $\<W,x\>^5$
in $(II)$, $(III)$ and $(IV)$ are respectively
  $$ 3 \frac{\sin \alpha}{r^2} +
(n-3) \frac{\sin \alpha}{r^2}$$
  $$  \frac{\sin \alpha}{r^2}- (n-1) \frac{\sin \alpha}{r^2}$$
  $$ -(n-1)\frac{ \sin \alpha}{r^2}- (n-1)^2 \frac{ \sin \alpha}{r^2},$$
from which we conclude the proof.
\end{proof}

\begin{lemm} \label{Alg}
   Let $P(x_1, x_2 ,x_3)$ be an irreducible polynomial with real
   coefficients. Assume
   the set of zeroes of $P$ is non empty and that
   $f \in  {\R}[x_1, x_2 , x_3]$ vanishes on it.
   Then $f$ is of the form $PQ$ with $Q \in {\R}[x_1, x_2 , x_3]$.
 \end{lemm}
\begin{proof}
 We embed $\R^3$ into $\C^3$. Assume by contradiction that $f$
 is not of the form $PQ$ with  $Q \in {\R}[x_1, x_2 , x_3]$.
 Then  also $f$ is of not the form $PQ$ with
 $Q \in {\C}[x_1, x_2 , x_3]$,
 so the set $Y:= \{x / f(x)=P(x)=0 \}=\{ f(x)=0 \} \cap \{ P(x)=0 \}$
 has complex codimension $2$, so
 real dimension $6-4=2$.
 Now our assumption is that $Y \cap \R^3$ contains $\{P=0\} \cap
 \R^3$, so in particular, $Y$ which has real dimension $2$ contains
 $\{P=0\} \cap
 \R^3$, also of real dimension $2$. We conclude that it is an
 irreducible component of $Y$, a contradiction because $\R^3$ does not
 contain any complex curve.
\end{proof}

\bigskip

 \noindent \emph{End of the proof of Theorem \ref{HS}.}\\
\emph{First case}: the three vectors $W$, $\dot{W}$ and
$\ddot{W}$ do not span $\R^n$ (it is in particular always the case
when $n > 3$).

In this case, there exist coordinates $(y_1, \ldots , y_n)$ on
$\R^n$ such that $f$ does not depend on $y_n$. In particular $\{y,
f(y)=0 \}$ contains some straight line
 $\{y_j= \mbox{Const},  1 \leq j \leq n-1 \}$.
 However by assumption $\{y, f(y)=0 \}$
 contains also the hypersphere $\S^{n-1}$. As it is an algebraic
 set, it is necessarily the whole $\R^n$, thus $f$ vanishes identically.
This implies that either $W$ vanishes (this is the centered
case), or that so does $f$ as a polynomial of \em independent \em
variables of $\R^n:$ the vectors $W$$,\dot{W}$ and $\ddot{W}$
might be non independent and in this case we should rewrite $f$ as
a polynomial of less variables. Anyway, we know from Lemma
\ref{ordre5} that the only term of order $5$ in $f$ is
 $\sin \alpha$ times a non negative constant, thus we deduce in
 any case that $\sin \alpha$ should vanish, which means that the
 curve $\gamma$ would be a line passing through the origin.
 Then we know that by Example 2 the immersion is totally geodesic.

\medskip

\noindent \emph{Second case}: 
$n=3$ and $\mbox{Span}(W,\dot{W},\ddot{W})=\R^3$.

Here we apply the algebraic
 Lemma \ref{Alg} with $P(x_1,x_2,x_3)=|x|^2-1$,
thus obtaining  that either $f$ is a multiple of $|x|^2-1$,
 or it vanishes. In the first case, there should be at least one
 term of degree $2$ in $\<\ddot{W},x\>$, which is not the case
 (\emph{cf} Lemma \ref{Wseconde}). So we deduce again that $f$ vanishes and
 the conclusion is the same as above.
 \end{proof}

\bigskip

\section{Centered Hamiltonian stationary $\sigma$-submanifolds}
\subsection{The differential system}
In this case, the induced metric is diagonal:
$g=\mbox{diag}(1,r^2, ... , r^2)$, and $\det g= r^{2(n-1)}$.
Moreover, $\beta$ depends only on $s:$ $\beta=\al + n \phi$, thus
$$ \Delta \beta= -\frac{1}{\sqrt{\det g}}\frac{\pa}{\pa s} 
          \left(\sqrt{\det g} \, \frac{\pa \beta}{\pa s} \right).$$
Hence the submanifold is Hamiltonian stationary if and only if
$$\sqrt{\det g} \, \frac{\pa \beta}{\pa s} =C,$$
which amounts to
$$r^{n-1}(n \dot{\phi}+ \dot{\alpha})=C,$$
for some non vanishing constant $C$
(the case of vanishing $C$ implies $\beta$ to be constant, so this
is the special Lagrangian case).
 We add that we could have used the
computations of the previous section as well with $W=0$ to find
the same equation.

Thus we are reduced to studying the following differential system
(using that $\dot{\phi}=\frac{\sin \alpha}{r}$):
$$ \left\{ \begin{array}{l}
 \dot{\alpha}=\frac{C}{r^{n-1}}-
      \frac{n \sin \alpha}{r}\\
     \dot{r}=\cos \alpha .  \end{array} \right. $$
There are fixed points
$(\bar{\alpha}=\pm  \frac{\pi}{2} \bmod 2\pi
  ,\bar{r}=\left(\frac{|C|}{n}\right)^{1/(n-2)})$,
corresponding to the case of $\gamma$ being a circle
centered at the origin. Up to a
sign change of parameter $s$, we may assume that $C$ is non
negative and then $\overline{\alpha}=\frac{\pi}{2} \bmod 2\pi$.
Moreover, the system admits a first integral: $E=2r^n \sin \alpha - Cr^2$,
so we can draw the phase portrait easily. As the system is
periodic in the variable $\alpha$, it is sufficient to study
integral lines in a strip of length $2 \pi$. Moreover integral lines
are symmetric with respect to the vertical lines ${\alpha=\pi/2 \bmod \pi}$.
  The energy of the fixed points is 
$E_0:=  \left( \frac{C}{n} \right)^{n/(n-2)}(2-n)$.
 There is a critical integral line of energy $E_0$
with bounded pieces connecting two fixed points and unbounded
pieces starting from (or ending to) a fixed point and having a 
branch asymptotic to a vertical line ${\alpha=k \pi, k \in \Z}$.
\begin{figure}[h]
\begin{center}
\input{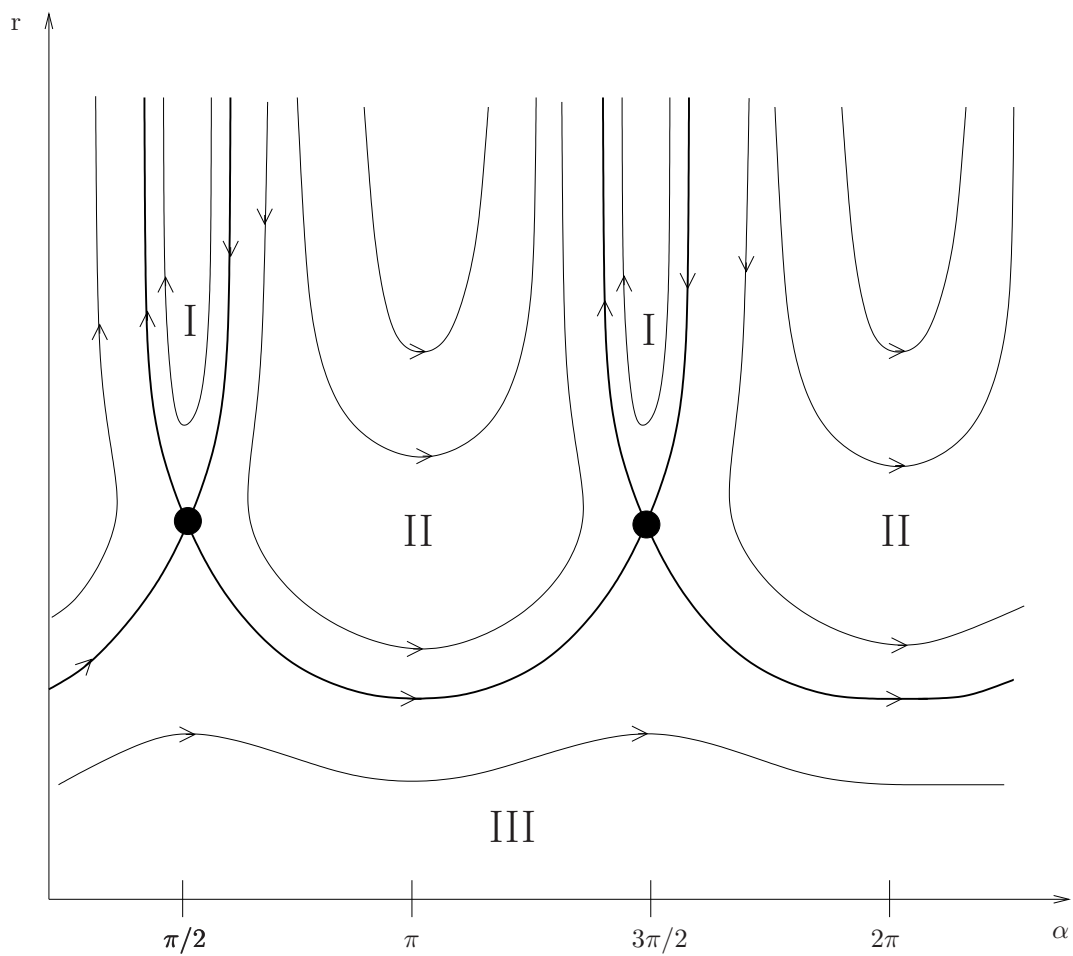}
\caption{Phase diagram}
\label{figure:phase}
\end{center}
\end{figure}
 We observe that the integral lines corresponding to a
non negative energy
  (resp. strictly less than $E_0$) are unbounded, we call them
\em type I \em (resp. \em type II \em) curves.
On the other hand, 
integral lines of energy $E \in (E_0,0)$ have two connected
components (up to periodicity), one of them being bounded in the
variable $r;$ this will be called in the following a \em type III \em curve.
As the unbounded component shares the same features as the curves of
non negative energy, it is also called a \em type I \em curve.

In order to have a better picture of the corresponding curves
$\gamma$ and of Hamiltonian stationary Lagrangian submanifolds
they generate, we shall discuss in the next two paragraphs 
 the inflection points of $\ga$ and
the
quantity 
$\Phi(E):=\int \dot{\phi}=\int \frac{\sin \alpha}{r}$
 that we call \em total variation of
phase. \em

\subsection{Inflection points of $\ga$}
They correspond to the vanishing of the curvature
$k=\dot{\theta}=\frac{C}{r^{n-1}}-\frac{(n-1)\sin \alpha}{r}$.
This implies the relation $ r=\big(\frac{C}{(n-1) \sin
\alpha}\big)^{1/(n-2)}$. 

We observe that in the case of dimension $3$, this is exactly the
equation of the integral curve of level $E=0:$ so in this
dimension there is  a solution which is of curvature
zero, i.e. a straight line, and the other ones don't have any
inflection point and are locally convex.

In higher dimension, let us compute the energy $E$ at points of
vanishing curvature. At such a point $(\alpha,r)$, we have
$r^{n-2}=\frac{C}{(n-1) \sin \alpha}$ so we
deduce that $E=C r^2 (\frac{2}{n-1}-1)$, which has range
$(-\infty,E_1)$ where $E_1$ is the energy level defined by 
 $E_1=2r_1^n - C r_1^2$ ($r_1$ is the least radius on the
curve of points $(\alpha,r)$ corresponding to vanishing curvature:
$r_1=\big(\frac{C}{(n-1) }\big)^{1/(n-2)}$).
As $r \mapsto E(r)=2r^n- C r^2$ is increasing when 
$r > \overline{r}$, $E_1$ belongs to the
interval $(E_0,0)$. We conclude that when $n >3$
every type II curve has two (symmetric)
 inflection points. This is also the case of some type I curves, 
while the remaining ones are locally convex.

\subsection{Study of the total variation of phase}
 We first compute $\Phi(E)$ for type I curves.
 Let $r_0$ be the minimal value
taken by $r$. We have the relation $E=2r_0^n-Cr_0^2$. We shall use the fact
that on the curve $\sin \alpha > 0;$ moreover, by symmetry we may restrict
ourselves to half of the curve, thus we may also 
assume that $\cos \alpha > 0$. Thus
$$\Phi=\int \frac{\dot{\phi}}{\dot{r}}dr=
2 \int_{r_0}^{\infty} \frac{\sin \alpha}{r \cos\alpha}dr
   =2 \int_{r_0}^{\infty} \frac{dr}{r}
 \left(\frac{1}{\sin^2 \alpha}-1\right)^{-1/2}$$
$$ =2 \int_{r_0}^{\infty} \frac{dr}{r}
   \left( \left(\frac{2r^n}{Cr^2+E} \right)^2 - 1 \right)^{-1/2}$$
   $$ =2 \int_{r_0}^{\infty} \frac{
   dr}{r}
   \left( \left(\frac{2r^n}{C(r^2-r_0^2)+2r_0^n} \right)^2 - 1 \right)^{-1/2}$$
Making the change of variable $r= x r_0$, we infer
$$\Phi=2 \int_{1}^{\infty} \frac{dx}{x}
   \left( \left(\frac{x^n}{\lambda(x^2-1)+1} \right)^2 - 1 \right)^{-1/2},$$
where $\lambda=\frac{C}{2r_0^{n-2}}$.
 We observe that this integral equals $\pi /n$ when $\lambda$ vanishes, is divergent for
 $\lambda=n/2$, \em i.e. \em for $E=E_0$ and is decreasing in the
 variable $r_0$, so $\Phi(E)$ has range
 $(\pi / n, + \infty)$.
Notice that the sign of $\dot{\phi}=\sin \alpha /r$ does not change, so 
that $\gamma=r e^{i \phi}$ is embedded whenever the total variation of 
phase is small enough, i.e. whenever $E$ is large enough.

On the other hand a
 unbounded piece of the
 integral curve of energy level $E_0$ is singular, with an
 infinite spiral branch asymptotic to the unit circle.

\bigskip

We now look at the case of type II integral curves.
Let $r_1$ be the value of
$r$ at $\alpha=0 \bmod \pi$,  so that
$ E = -C r_1^2$. 
We shall calculate separately the contributions of the curve when
the integrand is positive (resp. negative), corresponding to the part of
the curve lying in the region $\{\sin \al > 0\}$ 
(resp. $\{\sin \al < 0\}$).
$$\Phi_+:=2 \int_{r_1}^{\infty} \frac{\sin \alpha}{r \cos\alpha}dr
   =2 \int_{r_1}^{\infty} \frac{dr}{r}
 \left(\frac{1}{\sin^2 \alpha}-1\right)^{-1/2}$$
$$ =2 \int_{r_1}^{\infty} \frac{dr}{r}
   \left( \left(\frac{2r^n}{Cr^2+E} \right)^2 - 1 \right)^{-1/2}$$
   $$ =2 \int_{r_1}^{\infty} \frac{dr}{r}
   \left( \left(\frac{2r^n}{C(r^2-r_1^2)} \right)^2 - 1 \right)^{-1/2}$$
Making the change of variable $r= x r_1$, we obtain
$$\Phi_+=2 \int_{1}^{\infty} \frac{dx}{x}
   \left( \left(\frac{x^n}{\lambda(x^2-1)} \right)^2 - 1 \right)^{-1/2},$$
where $\lambda=\frac{C}{2r_1^{n-2}}$.

One calculates that $\Phi_+$ is increasing in the variable $E$
 and that $\lim_{E \to -\infty }= \pi / n$.

\bigskip

In the computation of $\Phi_-$, we may use $\al$ as a parameter:
$$\Phi_-:= \int_{\pi}^{2 \pi} \frac{\dot{\phi}}{ \dot{\al}}d\al
   = \int_{\pi}^{2 \pi}
  \frac{ \frac{\sin \al}{r} }{  \frac{C}{r^{n-1}} - \frac{n \sin \al}{r} }
                      d \al$$
   $$ = \int_{\pi}^{2 \pi} 
    \frac{r^{n-2} \sin \al}{C- n r^{n-2}\sin \alpha} d \al$$
It is an easy computation to show that 
$- \pi /n < \Phi_- < 0$ and that
$\lim_{E \to -\infty} \Phi_-= -\pi / n$.

Next we show that $\gamma$ has always
a self-intersection. Let $s_0$ be the parameter value corresponding
to the point of the curve of least radius (in particular $\al(s_0)= 3 \pi /2 \bmod 2\pi$).
The  fact that $\Phi_+ > |\Phi_-|$ and the intermediate value theorem 
imply that there exists $s_1$ such that
$\phi(s_1)=\phi(s_0);$ 
moreover by the symmetry of the phase portrait, there exists 
$s_2 \neq s_1$ with the same property, the corresponding points
in the phase portrait satisfying $r(s_1)=r(s_2)$ and  
  $ \frac{1}{2} (\al(s_1)+\al(s_2))= 3 \pi / 2 \bmod 2\pi$. Thus $\ga(s_1)=\ga(s_2)$.

\bigskip

We end this section by looking at the type III curves. 
As $\alpha$ is always increasing on the curve, we shall use it as
a parameter and consider the piece of  curve 
$\ga([\pi/2, 2 \pi + \pi / 2])$.

 As $\phi$ is decreasing on $\alpha \in [\pi,2 \pi]$ and
 increasing elsewhere, we have:

$$\Phi_+:= \int_{\pi/2}^{\pi} \dot{\phi}
 + \int_{2 \pi}^{2\pi+ \pi / 2}  \dot{\phi} > 0 $$ and
$$\Phi_-:= \int_{\pi}^{2 \pi}   \dot{\phi} <0. $$

It is easy to show (the calculations are left
to the Reader) that $\Phi_+ > |\Phi_-|$ and that
 $\lim_{E \to 0} \Phi_+ = \lim_{E \to 0} \Phi_- = 0$.
Moreover, $\lim_{E \to E_0} \Phi_+ = + \infty$ and
$\lim_{E \to E_0} \Phi_-$ is finite. This implies that
the range of the total variation of phase for type III curves
is $(0,\infty)$. In particular, to the limiting case
$E=E_0$ corresponds a bounded, complete curve spiraling 
asymptotically to the unit circle.

With a similar argument as above, we show that here again $\ga$ has
always self-intersections.
We conclude that for a type III curve  such
that $\Phi(E) \in 2 \pi \mathbb{Q}$, the corresponding curve
$\gamma$ is bounded, closed and non-embedded.

%

\subsection{Conclusion (Corollary 1)}
 We are now in position to describe the whole family of
 Hamiltonian stationary $\si$-submanifolds:
\begin{coro} Any Hamiltonian stationary Lagrangian submanifold which
is foliated by $(n-1)$-spheres is locally congruent to one of the following:

\begin{itemize}
\item The standard embedding
$$ \begin{array}{ccc}\S^1 \times \S^{n-1}  & \to & \C^n  \\
(e^{is} , x) & \mapsto &  e^{is}x .
\end{array}$$

\item A singular, bounded  immersion of $\R \times \S^{n-1}$
``spiraloid" asymptotic to $\S^1 \times \S^{n-1}$.

\item A singular, unbounded  ``spiraloid" with a
 smooth end and asymptotic to $\S^1 \times \S^{n-1}$.

\item A family of smooth ``catenoid-type" immersions of $\R  \times \S^{n-1}$,
some of them being embedded.
 In dimension $3$, one of them takes the particular following form:
 $$ \begin{array}{ccc}\R \times \S^{2}  & \to & \C^3  \\
(s , x) & \mapsto &  (C+is)x,
\end{array}$$
where $C$ is some real constant.

\item A family of non-standard smooth immersions of $\S^1 \times \S^{n-1}$.
   They always have self-intersections.

\end{itemize}
\end{coro}

\section*{Appendix: computation of $\Delta \beta$ (proof of Lemma \ref{computation})}
 Here we are going to use the same notations as in Section 2 and 4.
We first compute the following
$$ \<\ell_{ss},\ell_{s}\>=  \<W, \dot{W}\>
  + \left(  \frac{ \sin \alpha }{r} -
        k \sin \alpha \right) \<W,x\>+ \cos \alpha \<\dot{W},x\>$$
$$ \<\ell_{ss},\ell_{*}v_{j}\>=r \<\dot{W}, v_{j}\>$$
$$ \<v_{j}(\ell_{s}),\ell_{s}\>= \cos \alpha \<W,v_{j}\>$$
$$ \<v_{j}(\ell_{s}),\ell_{*}v_{k}\>=\delta_{jk} r \cos \alpha $$
$$ \<\ell_{v_{j}v_{k}},\ell_{s}\>=-\delta_{jk}r(\cos \alpha + \<W,x\>)$$
   $$\<\ell_{v_{j}v_{k}},\ell_{*}v_{l}\>=0.$$

 We shall make use of the following formula:
$$\Delta \beta = \mbox{div} J n\vec{H}=
    \< \nabla_{e_{1}} J n\vec{H}, e_{1} \> +
    \sum_{j=2}^{n} \< \nabla_{e_{j}} J n\vec{H}, e_{j} \>$$
Writing $Jn\vec{H}= a e_{1}+\sum_{k} a_{k} e_{k}$, we have
$$ \nabla_{e_{j}} J n\vec{H}=
       e_{j}(a)e_{1}+a \nabla_{e_{j}}e_{1} +
     \sum_{k} e_{j}(a_{k})e_{k}+ \sum_{k} a_{k} \nabla_{e_{j}}e_{k}, $$
 and an analogous expression for $\nabla_{e_1} J n\vec{H}$.

 Then we obtain
\begin{eqnarray*}
\Delta \beta &=& e_{1}(a) + a \< \nabla_{e_{1}} e_{1}, e_{1} \>
      + \sum_{k} a_k \< \nabla_{e_1} e_k, e_1 \>
\\
&& + \sum_j \left[ a \< \nabla_{e_j} e_1, e_j \> + e_j(a_j)
     +  \sum_k  a_k \< \nabla_{e_j} e_k, e_j \>  \right]. 
\end{eqnarray*}

We give the expressions of the coefficients of $J\vec{H}:$

$$ a_{j}= A^2 \frac{\sin \alpha \<W, v_j\>}{r} $$
     and
     $$a =- A^3 B  -A (n-1) \frac{\sin \alpha}{r},$$
  where
 $$ A=\left(1+ 2  \cos \alpha \<W,x\>+\<W,x\>^2 \right)^{-1/2}$$
  and
  $$B=k+ \left(k \cos \alpha  +
  \frac{\sin \alpha \cos \alpha}{r} \right) \<W,x\>+
  \sin \alpha \< \dot{W},x\> +
 \frac{\sin \alpha}{r}\<W,x\>^2 .$$

 We start with some easy computations:

 \begin{lemm}

   \begin{description}

    \item[a)] $\dot{A}=-A^3 \left(
    \cos \alpha \< \dot{W},x\>
    - \dot{\alpha} \sin \alpha \<W,x\>
    + \<W,x\> \< \dot{W},x\> \right)$,
    \item[b)] $\<\pa_s,e_1 \> = A^{-1}, $
    \item[c)] $v_j(A)=-A^{3}\Big(\cos \alpha + \< W,x\> \Big)\<W, v_j \>$,
     \item[d)] $v_j(B)= \Big(
  k \cos \alpha + \frac{\sin \alpha \cos \alpha}{r}
   +  \frac{2 \sin \alpha}{r}\<W,x\> \Big) \<W,v_j\>
  +\sin \alpha \< \dot{W},v_j\>, $
     \item[e)] $ \nabla_{v_j}v_{k}=
        -\delta_{jk}Ar(\cos \alpha + \<W,x\>) e_{1}, $
    \item[f)] $ \nabla_{v_j} \pa_s = \cos \alpha \; e_j$.
    \end{description}
\end{lemm}

\begin{proof}
\textbf{a),b),c)} and \textbf{d)} are easy and left to the reader.
\begin{description}
\item[e)]
   We write $\nabla_{v_{j}} v_{k}=b_{1} \pa_{s}+\sum_m b_{m} v_{m}$.
    In particular, we have
  $$ 0=\<\nabla_{v_{j}} v_{k},v_{l}\>=
       b_1 \<\pa_s, v_l|> + \sum_m b_m \<v_m, v_l |> = 
           b_{1} g_{1l}+ \sum_m b_{m} g_{ml}= b_{1}g_{1l}+ b_{l}r^{2},$$
  so we obtain the relation $ b_{l}= - \frac{\<W,v_{l}\>}{r}b_{1}$.
  
  When $j \neq k$, it yields
  $$ 0 = \<\nabla_{v_{j}} v_{k}, \pa_{s}\>=
    b_{1} \left(g_{11}-\sum_l  \frac{\<W,v_{l} \>}{r} g_{1l} \right),$$
      so all the coefficients $b_{1}, b_{l}$ vanish, and we obtain
  \textbf{e)}.

  If $j=k$, we have
   $$ -r( \cos \alpha + \<W,x\>)
   = \<\nabla_{v_j} v_{j}, \pa_{s}\>=
       b_{1} \left(g_{11}-\sum_l  \frac{\<W,v_{l} \>}{r} g_{1l} \right),$$
 so $b_{1}=-A^{2}r( \cos \alpha + \<W,x\>)$.
 We deduce that
 $$ \nabla_{v_j} v_{j}=b_{1} \left(\pa_{s}-
        \sum_l \frac{\<W,v_l\>}{r}v_l \right)
        =b_{1} \frac{e_1}{A}
        =-Ar( \cos \alpha + \<W,x\>)e_{1},$$
 which implies  \textbf{e)}.

\item[f)]
   We write $\nabla_{v_j} \pa_s=b_{1} \pa_s+\sum_l b_l v_l$.
 Then we have
 $$\delta_{jk} r \cos \alpha =
   \<\nabla_{v_j} \pa_s, v_k \> =b_1 g_{1k}+\sum_l b_l g_{kl}.$$

  When $j \neq k$, we deduce that
    $$ b_k= -\frac{g_{1k}}{g_{kk}} b_1=-\frac{\<W,v_k\>}{r} b_1$$

 When $j = k$, we obtain
 $$r \cos \alpha = b_1 r \<W,v_j\> + r^{2} b_j,$$
 so 
 $$ b_j= r^{-2} ( r \cos \alpha - b_1 r \<W,v_j\>)=
               r^{-1} (  \cos \alpha - b_1  \<W,v_j\>)  .$$

 From these relations we can write
$$  \cos \alpha  \<W,v_k\> = \<\nabla_{v_k} \pa_s, \pa_s \>
      =  b_1 g_{11} +\sum_j b_j g_{1j} $$
$$= b_1 \left( g_{11} - \sum_j r \<W,v_j\> \frac{\<W,v_j\>}{r}\right)
       +  r \<W,v_k\>  r^{-1} \cos \alpha
     = b_1 A^{-2} +  \cos \alpha \<W,v_k\>$$
 We deduce that $b_1$ vanishes and also $b_l$ for $l \neq k$,
 and finally $b_k= \frac{\cos \alpha}{r}$, so we obtain \textbf{f)}.
\end{description}
\end{proof}
\bigskip

 We now compute the different terms of $\Delta \beta$:

\begin{lemm}
   \begin{description}

    \item[a)]  $$ e_1(a)=A \dot{a} - \sum_j \frac{A \<W,v_j\>}{r}
    v_j(a)$$
     with

    $$ A \dot{a} = 3 A^6 B  \Big(\cos \alpha \< \dot{W},x\>
    - \dot{\alpha} \sin \alpha \<W,x\>
    + \<W,x\> \< \dot{W},x\> \Big) $$
    $$-
    A^4 \left[\dot{B} - (n-1) \frac{ \sin \alpha}{r}
  \Big(\cos \alpha \< \dot{W},x\>
    - \dot{\alpha} \sin \alpha \<W,x\>
    + \<W,x\> \< \dot{W},x\> \Big) \right]$$
$$ - A^2 (n-1) \pa_{s} \left(\frac{\sin \alpha}{r} \right)$$

and

    $$ \sum_j \frac{A \<W,v_j\>}{r} v_j(a)=
   3 A^6 B \frac{ \cos \alpha +\<W,x\>}{r}(|W|^2-\<W,x\>^2)  $$
     $$ - \frac{A^4}{r} \left[\left(
       k \cos \alpha
  - (n-2) \frac{\sin\alpha\cos\alpha}{r}
  - (n-3) \frac{ \sin \alpha}{r}  \<W,x\>
  \right)  (|W|^2-\<W,x\>^2) \right.$$
 $$\left.   +\sin \alpha \Big(\<\dot{W},W\>-\<\dot{W},x\> \<W,x \> \Big)  \right],$$

    \item[b)] $$\< \nabla_{e_1} e_1, e_1 \>=0, $$

     \item[c)] 
\[ \begin{array}{rrl}
\sum_{k} a_k \< \nabla_{e_1} e_k, e_1 \>+\sum_j  e_j(a_j) = 
- \frac{\sin \alpha}{r^2} \Big( 
& A^{4} (\cos \alpha + \<W,x\> ) ( |W|^{2}- \<W,x\>^{2})
\\
&+  A^{2} (n-1) \<W,x\> & \Big),
\end{array} \]

    \item[d)]  $$\sum_j \<\nabla_{e_j} e_1, e_j \>=
    \frac{A(n-1)}{r}\left( \cos \alpha +\<W,x\>\right),$$

    \item[e)] $$\<\nabla_{e_j} e_k, e_j \>=0, \quad  \forall j,k.$$
   \end{description}
\end{lemm}
\begin{proof}
\begin{description}
\item[a)]
 We compute
 \begin{eqnarray*}
\dot{a} &=& 
- 3 A^2 \dot{A} B
  - A^{3} \dot{B} - \dot{A} (n-1) \frac{\sin \alpha}{r}
   - A \pa_{s} \left((n-1)\frac{\sin \alpha}{r} \right)
\\
&=& A^3\left(3A^2 B + (n-1) \frac{ \sin \alpha}{r} \right)
  \Big(\cos \alpha \< \dot{W},x\>
    - \dot{\alpha} \sin \alpha \<W,x\>
    + \<W,x\> \< \dot{W},x\> \Big)
\\
&& - A^3 \dot{B} - A (n-1) \pa_{s} \left(\frac{ \sin \alpha}{r} \right)
\end{eqnarray*}

from which we deduce the first equality. Next we have
\begin{eqnarray*}
v_j(a) &=& 
-3 A^2 v_j(A) B  - A^{3} v_j(B) - v_j(A) (n-1)\frac{\sin \alpha}{r}
\\
&=& A^{3} \left(3 A^2B + (n-1) \frac{\sin\alpha}{r}\right)
                           \Big( \cos \alpha + \<W,x \> \Big) \<W,v_j\> 
\\
&& - A^{3} \left[\left(
  k \cos \alpha + \frac{\sin \alpha \cos \alpha}{r}
  +\frac{2 \sin \alpha}{r}\<W,x\> \right) \<W,v_j\>
   +\sin \alpha \< \dot{W},v_j\>  \right]
\\
&=&  3 A^{5}( \cos \alpha + \<W,x \>) B \<W,v_j\>
\\
&& - A^{3} \left[\left(
   k \cos \alpha + \frac{\sin \alpha \cos \alpha}{r}
  -(n-1) \frac{\sin \alpha \cos \alpha}{r} \right.\right. 
\\
&& \left. \left. +\left(\frac{2 \sin \alpha}{r} - (n-1)\frac{ \sin \alpha}{r} 
\right) \<W,x\>  \right) \<W,v_j\> +\sin \alpha \< \dot{W},v_j\>  \right]
\\
&=&  3 A^{5}( \cos \alpha + \<W,x \>) B \<W,v_j\>
\\
&& -A^{3} \Bigg[\left(
   k \cos \alpha -(n-2) \frac{\sin \alpha \cos \alpha}{r}
  - (n-3)\frac{ \sin \alpha}{r}  \<W,x\>
  \right) \<W,v_j\>
\\
&& +\sin \alpha \< \dot{W},v_j\>  \Bigg] .
\end{eqnarray*}
Thus
\begin{multline*}
A \sum \frac{v_j(a) \<W,v_j\>}{r} =
     3 A^6 B \frac{ \cos \alpha +\<W,x\>}{r} \sum \<W,v_j\>^2 
\\
 -\frac{A^4}{r} \Bigg[\left(
       k \cos \alpha
  - (n-2) \frac{\sin\alpha\cos\alpha}{r}
  - (n-3) \frac{ \sin \alpha}{r}  \<W,x\>
  \right) \sum \<W,v_j\>^2
\\
   +\sin \alpha \sum \<W,v_j\> \< \dot{W},v_j\>  \Bigg]
\end{multline*}

 \item[b)]
 We first compute
\begin{eqnarray*}
\nabla_{e_1} e_1 &=& 
A \Bigg[
  \nabla_{\pa_s} \left(A\pa_s -\frac{A}{r} \sum \<W,v_j\> v_j \right)
\\
&&    - \sum \frac{\<W,v_j\>}{r} \nabla_{v_j}
       \left( A \pa_s - \frac{A}{r} \sum \<W,v_k\> v_k \right)
       \Bigg]
\\
&=&   A \Bigg[
      A \nabla_{\pa_s} \pa_s + \dot{A} \pa_s
      - \sum \pa_s \left( \frac{A \<W,v_j\>}{r}\right) v_j
   - \frac{A}{r} \sum \<W,v_j\> \nabla_{\pa_s} v_j
\\
&& - \sum \frac{ \<W,v_j\>}{r} A \nabla_{v_j} \pa_s  
  - \sum \frac{ \<W,v_j\> v_j(A)} {r} \pa_s
\\
&&     + \frac{1}{r^2} \sum_{j,k} \<W,v_j\> \Big(
   v_j(A)  \<W,v_k\> v_k + A( -\delta_{jk}  \<W,x\>v_k
       +A  \<W,v_k\> \nabla_{v_j} v_k \Big) \Bigg].
 \end{eqnarray*}
 Making the scalar product with $e_1$, many terms vanish and we obtain
\begin{eqnarray*}
\<\nabla_{e_1} e_1,e_1\> &=& 
A^2 \<\nabla_{\pa_s} \pa_s,e_1\>     + A \dot{A} \<\pa_s,e_1 \>
  - \frac{A}{r} \left( \sum \<W,v_j\>v_j(A)\right)\<\pa_s,e_1 \>
\\
&& -\frac{A^{3}}{r} \left(\cos \alpha  +  \<W,x\> \right)\sum \<W,v_j\>^{2}
\\
&=& A^3 \left( \dot{\alpha} \sin \alpha \<W,x\> 
+ \frac{\sin^2 \alpha \<W,x\>}{r} -      k \sin
\alpha  \<W,x\> \right) 
\\
&& +\frac{1}{r} \left( \sum \<W,v_j\>^2 A^3(\cos \alpha + \<W,x\>) \right)
\\
&& -\frac{A^{3}}{r} \left(\cos \alpha  +  \<W,x\> \right)\sum \<W,v_j\>^{2}.
\end{eqnarray*}
Using the fact that
  $\dot{\alpha}=\dot{\theta}-\dot{\phi}=k - \frac{\sin \alpha}{r}$
we see that the first term in the last sum vanishes, and so does the last two
ones together, so we obtain  \textbf{b)}.

\item[c)]
  It will be a consequence of
  the fact that $|W|^2= \<W,x\>^2 + \sum \<W,v_j\>^{2}$
  and of the two following equalities:
\begin{multline} \sum e_j(a_j) = -\frac{\sin \alpha}{r^2}
    \bigg(  2A^{4} (\<W,x\>+ \cos \alpha \>) \sum \<W,v_j\>^{2}
\\
    + (n-1) A^{2} \<W,x\> \bigg),    \tag{1}
\end{multline}
\[
\sum a_k \<\nabla_{e_1} e_k, e_1 \> =
       \frac{ A^4}{r^2} \sin \alpha
      ( \cos \alpha +  \<W,x\>) \sum \<W,v_k\>^{2}. \tag{2}
\]
  To prove $(1)$ and $(2)$ we proceed as follows:
 $$ e_j(a_j)=\frac{1}{r} v_j
    \left(A^{2} \frac{\sin \alpha \<W,v_j\>}{r} \right)
     = \frac{\sin \alpha}{r^2}
    \left( 2A v_j(A) \<W,v_j\>- A^{2} \<W,x\> \right)  $$
 $$=-\frac{\sin \alpha}{r^2}
    \left(  2A^{4} (\<W,x\>+ \cos \alpha \>) \<W,v_j\>^{2}
    + A^{2} \<W,x\> \right),  $$
from which we obtain $(1)$. Then we compute
\begin{eqnarray*}
\nabla_{e_1} e_k &=&
A \left( \nabla_{\pa_s} e_k -
     \sum_{l} \frac{\<W,v_l\>}{r}  \nabla_{v_l}e_k\right)
\\
&=& A \left(\frac{1}{r} \nabla_{\pa_s} v_k +
             \pa_{s}(r^{-1}) v_k+
     \sum_{l} \frac{\<W,v_l\>}{r^2}  \nabla_{v_l}v_k\right)
\\
&=& A \left(\frac{ \cos \alpha}{r} e_k -
        \frac{\cos \alpha}{r^2}  v_k +
     \sum_l \frac{\<W,v_l\>}{r^2}  \nabla_{v_l}v_k\right). 
\end{eqnarray*}
  The first two terms of the last expression
vanish and thanks to the
  \textbf{e)} of Lemma 6, we deduce
  $$ \<\nabla_{e_1} e_k, e_{1} \> =
      - \frac{ A^2\<W,v_j \>}{r} ( \cos \alpha +  \<W,x\>),$$
  which implies $(2)$.

\item[d)]
We have
\begin{eqnarray*} 
\nabla_{e_j} e_{1} &=&
\frac{1}{r} \nabla_{v_j} \left(
A \pa_{s} - \sum_k \frac{A\<W, v_{k}\>}{r} v_{k} \right)
\\
&=& \frac{1}{r} \left[ v_{j}(A) \pa_{s} + A \nabla_{v_j}\pa_{s}
    - \sum_k v_{j} \left( \frac{A \<W,v_{k} \>}{r} \right) v_{k}
   - \sum_k \frac{A \<W,v_{k} \>}{r}  \nabla_{v_{j}} v_{k}
   \right]
\\
&=&\frac{1}{r} v_{j}(A) \pa_s + A \frac{\cos \alpha}{r}e_j
  - \frac{1}{r} \sum_k  v_j\left( \frac{A \<W,v_{k} \>}{r}\right) v_k
\\
&&     + \frac{A \<W,v_j \>}{r^2}r(\cos \alpha+  \<W,x\>)e_{1}
\end{eqnarray*}
so
\begin{eqnarray*}
\< \nabla_{e_j} e_{1},e_{j} \>  &=&
   \frac{1}{r} v_j(A) \<\pa_s,e_{j}\> + A \cos \alpha
    -  v_{j} \left( \frac{A \<W,v_{j} \>}{r} \right)
\\
&=& \frac{1}{r}v_{j}(A) \<W,v_j\> +   A \cos \alpha
       -v_{j}(A) \frac{\<W,v_{j}\>}{r} + \frac{A \<W,x\>}{r}
\end{eqnarray*}
The first and the third terms cancel out, and we obtain \textbf{d)}.

\item[e)]
 This is an obvious consequence of part \textbf{e)} the Lemma 6.
\end{description}
\nopagebreak
\end{proof}
\medskip

Summing these computations, and using the formula for $\Delta \beta$
given at the beginning of this Appendix, we conclude that $f:=A^{-6} \Delta
\beta=(I)+(II)+(III)+(IV)$ with
 \begin{multline*}
 (I):=A^{-5} \dot{a} = 3 B  \Big(\cos \alpha \< \dot{W},x\>
    - \dot{\alpha} \sin \alpha \<W,x\>
    + \<W,x\> \< \dot{W},x\> \Big)
\\
- A^{-2} \left[\dot{B} - (n-1) \frac{ \sin \alpha}{r}
  \Big(\cos \alpha \< \dot{W},x\>
    - \dot{\alpha} \sin \alpha \<W,x\>
    + \<W,x\> \< \dot{W},x\> \Big) \right]
\\
- A^{-4} (n-1) \pa_{s} \left(\frac{\sin \alpha}{r} \right)
\end{multline*}
 \begin{eqnarray*}
 (II) &:=& -A^{-6} \sum \frac{\<W,v_j\> v_j(a)}{r}
\\
&=& -3  B \frac{ \cos \alpha +\<W,x\>}{r}(|W|^2-\<W,x\>^2)
\\
&& + \frac{A^{-2}}{r} \left[\left( k \cos \alpha
  - (n-2) \frac{\sin\alpha\cos\alpha}{r}
  - (n-3) \frac{ \sin \alpha}{r}  \<W,x\>
  \right)  (|W|^2-\<W,x\>^2) \right.
\\
&& \left.   +\sin \alpha \Big(\<\dot{W},W\>-\<\dot{W},x\> \<W,x \> \Big)  \right]
\end{eqnarray*}
\begin{eqnarray*}
(III) &:=& A^{-6} \left( \sum_{k} a_k \< \nabla_{e_1} e_k, e_1 \>
   + \sum_j  e_j(a_j) \right)
\\
&=& -A^{-2} \frac{\sin \alpha}{r^2}
     (\cos \alpha + \<W,x\> ) \Big( |W|^{2}- \<W,x\>^{2} \Big) -
     A^{-4}(n-1) \frac{ \sin \alpha}{r^2} \<W,x\>
\end{eqnarray*}
\begin{eqnarray*}
(IV) &:=&  A^{-6} a \sum \<\nabla_{e_j} e_1, e_j \>
\\
&=& - A^{-2}B \frac{n-1}{r} (\cos \alpha + \<W,x\>)-
   A^{-4} \frac{(n-1)^{2} \sin \alpha}{r^2}  (\cos \alpha + \<W,x\>).
\end{eqnarray*}

\bigskip

Henri Anciaux,

IMPA, Estrada Dona Castorina, 110 

22460--320  Rio de Janeiro, Brasil

email address: {\tt henri@impa.br}

\bigskip

Ildefonso Castro

Departamento de Matem\'aticas

Universidad de Ja\'en,

23071 Ja\'en, Spain

email address: {\tt icastro@ujaen.es}

\bigskip

Pascal Romon

Universit\'e de Marne-la-Vall\'ee

5, bd Descartes, Champs-sur-Marne,

77454 Marne-la-Vall\'ee cedex 2, France

email address: \texttt{romon@univ-mlv.fr}

\end{document}